\documentclass[12pt, a4paper]{amsart}
\usepackage{amscd,amssymb,eucal,eufrak,mathrsfs,amsmath}

\usepackage{hyperref}

\input xy
\xyoption{all}

\newtheorem{thm}{Theorem}[section]

\newtheorem{lemma}[thm]{Lemma}
\newtheorem{prop}[thm]{Proposition}

\theoremstyle{definition}
\newtheorem{defin}[thm]{Definition}

\theoremstyle{remark}

\newtheorem{exam}[thm]{Example}

\numberwithin{equation}{section} 

\def\N{{\Bbb N}}

\newcommand{\hsk}{\hskip 0.3cm}

\newcommand{\sk}{\vskip 0.4cm}

\newcommand{\ben}{\begin{enumerate}}

\newcommand{\een}{\end{enumerate}}
\newcommand{\bit}{\begin{itemize}}

\newcommand{\eit}{\end{itemize}}

\def\eps{\varepsilon}

\def\QED{\nobreak\quad\ifmmode\roman{Q.E.D.}\else{\rm Q.E.D.}\fi}

\begin{document}

\title[]
{A note on sensitivity of semigroup actions}
\author[]{Eduard Kontorovich}
\address{Department of Mathematics,
Bar-Ilan University, 52900 Ramat-Gan, Israel}
\email{bizia@walla.co.il}
\urladdr{http://www.math.biu.ac.il/$^\sim$kantore}

\author[]{Michael Megrelishvili}
\address{Department of Mathematics,
Bar-Ilan University, 52900 Ramat-Gan, Israel}
\email{megereli@math.biu.ac.il}
\urladdr{http://www.math.biu.ac.il/$^\sim$megereli}

\date{July 18, 2006}


\begin{abstract}
It is well known that for a transitive dynamical system $(X,f)$
sensitivity to initial conditions follows from the assumption that
the periodic points are dense. This was done by several authors:
Banks, Brooks, Cairns, Davis and Stacey \cite{Br}, Silverman
\cite{Sil} and Glasner and Weiss \cite{GW}. In the latter article
Glasner and Weiss established a stronger result (for compact
metric systems) which implies that a transitive non-minimal
compact metric system $(X,f)$ with dense set of almost periodic
points is sensitive.
 This is true also for group actions as was
 proved in the book of Glasner \cite{Gl}.

 Our aim is to generalize these results in the frame
 of a unified approach for a wide class
of topological semigroup actions including one-parameter semigroup
actions on Polish spaces.
\end{abstract}

\maketitle

\section{Introduction}

First we recall some well known closely related results regarding
sensitivity of dynamical systems.

\begin{thm} \label{facts}
\ben
\item \emph{(Banks, Brooks, Cairns, Davis and Stacey \cite{Br};
Silverman \cite{Sil}\footnote{under different but very close
assumptions})} Let $X$ be an infinite metric space and $f: X \to
X$ be continuous. If $f$ is topologically transitive and has dense
periodic points then $f$ has sensitive dependence on initial
conditions.
\item
\emph{(Glasner and Weiss \cite[Theorem 1.3]{GW}); see also Akin,
Auslander and Berg \cite{Ak})} Let $X$ be a compact metric space
and the system $(X,f)$ is an \emph{M-system} and
not minimal.
Then $(X,f)$ is sensitive.
\item
\emph{(Glasner \cite[Theorem 1.41]{Gl})} Let $X$ be a compact
metric space. An almost equicontinuous M-system $(G,X)$, where $G$
is a group, is minimal and equicontinuous. Thus M-system which is
not minimal equicontinuous is sensitive. \een
\end{thm}

\emph{Topological transitivity} of $f: X \to X$ as usual, means
that for every pair $U$ and $V$ of nonempty subsets of $X$ there
exists $n > 0$ such that $f^n(U) \cap V$ is nonempty. Analogously
can be defined general semigroup action version (see Definition
\ref{Top trans}.1).

 If $X$ is a compact metric space then (2) easily covers (1).
 In order to explain this recall that \emph{M-system} means that the
set of almost periodic points is dense in $X$ (\emph{Bronstein
condition}) and, in addition, the system is topologically
transitive. A very particular case of Bronstein condition is that
$X$ has dense periodic points (the so-called \emph{P-systems}). If
now $X$ is infinite then it cannot be minimal.

 Our aim is to provide a unified and generalized approach.
We show that (2) and (3) remain true for a large class of
\emph{C-semigroups} (which contains: cascades, topological groups
and one-parameter semigroups) and M-systems (see Definitions
\ref{d:almost group} and \ref{d:minimality}). Our approach allows
us also to drop the compactness assumption of $X$ dealing with
Polish phase spaces. A topological space is \emph{Polish} means
that it admits a separable complete metric.

We cover also (1) in the case of Polish phase spaces.

Here we formulate one of the main results (Theorem \ref{main}) of
the present article.

\sk

\noindent \textbf{Main result:}  \emph{Let $(S,X)$ be a dynamical
system where $X$ is a Polish space and $S$ is a C-semigroup. If
$X$ is an $M$-system which is not minimal or not equicontinuous.
Then $X$ is sensitive.}

\sk

\noindent \textbf{Acknowledgment.} We thank E. Akin for helpful
suggestions and S. Kolyada and L. Snoha for sending us their
article \cite{Ko}.

\sk
\section{preliminaries}
 \sk

 A \emph{dynamical
system} in the present article is a triple $(S,X,\pi)$, where $S$
is a topological semigroup, $X$ at least is a Hausdorff space and
$$\pi: S \times X \to X, \ \ (s,x) \mapsto sx$$ is a continuous
action on $X$. Thus, $s_1(s_2x)=(s_1s_2)x$ holds for every triple
$(s_1,s_2,x)$ in $S\times S \times X$. Sometimes we write the
dynamical system as a pair
$(S,X)$ or even as $X$, when $S$ is understood. 
The \emph{orbit} of $x$ is the set $Sx:=\{sx: s\in S\}$.
By $\overline{A}$ we will denote the closure of a subset $A
\subset X$. If
 $(S,X)$ is a system and $Y$ a closed $S$-invariant subset, then we
 say that $(S,Y)$, the restricted action, is a \emph{subsystem} of
 $(S,X)$.  For $U \subset X$ and $s \in S$ denote
  $$s^{-1}U: = \{x\in X:sx\in U \}.$$

If $S=\{f^n\}_{n\in \N}$ (with $\N:=\{1,2, \cdots\}$) and $f: X
\to X$ is a continuous function, then the classical dynamical
system $(S,X)$ is called a \emph{cascade}. Notation: $(X,f)$.


 \begin{defin} \label{d:almost group}
Let $S$ be a topological semigroup.
 \item $(1)$  We say that $S$ is a (left)
\emph{$F$-semigroup} if for every $s_0 \in S$ the subset $S
\setminus Ss_0$ is finite. \item  $(2)$  We say that $S$ is a
\emph{$C$-semigroup} if $S \setminus Ss_0$ is relatively compact
(that is, its closure is compact in $S$).
\end{defin}

\begin{exam} \label{e:C}
\ben
\item Standard one-parameter semigroup $S:=([0, \infty), +)$ is a C-semigroup.
\item Every cyclic "positive" semigroup $M:=\{s^{n}:n\in\N\}$ is an
F-semigroup. In particular, for every cascade $(X,f)$ the
corresponding semigroup $S=\{f^n\}_{n\in \N}$ is an F-semigroup
(and hence also a C-semigroup).
\item Every topological group is of course an F-semigroup.
\item Every compact semigroup is a C-semigroup.
 \een
\end{exam}

\begin{defin}
Let $(S,X)$ be a dynamical system where $(X,d)$ is a metric space.
\ben
\item A subset $A$ of $S$ acts \emph{equicontinuously}
at $x_0 \in X$ if for every $\epsilon>0$ there exists $\delta>0$
such that $d(x_0,x)< \delta$ implies $d(ax_0,ax)<\epsilon$ for
every $a\in A.$
\item
A point $x_0 \in X$ is called an \emph{equicontinuity point}
(notation: $x_0 \in Eq(X)$) if $A:=S$ acts equicontinuously at
$x_0$. 
If $Eq(X)=X$ then $(S,X)$ is \emph{equicontinuous}.
\item
$(S,X)$ is called \emph{almost equicontinuous} (see \cite{Ak, Gl})
if the subset $Eq(X)$ of equicontinuity points is a dense subset
of $X$. \een
\end{defin}

\begin{lemma} \label{com-of-S}
Let $(S,X)$ be a dynamical system where $(X,d)$ is a metric space.
Let $A \subset S$ be a relatively compact subset. Then $A$ acts
equicontinuously on $(X,d)$.
%
\end{lemma}


\sk

\section{Transitivity conditions of semigroup actions}

\sk

 \begin{defin} \label{Top trans} The dynamical system $(S,X)$ is called:
 \begin{enumerate}
   \item \emph{topologically
 transitive} (in short: TT) if for every pair $(U,V)$ of
 non-empty open sets $U,V$ in $X$ there exists $s\in S$ with
 $U \cap sV \neq \varnothing$. Since
 $s(s^{-1}U \cap V)=U \cap sV$,
 it is equivalent to say that $s^{-1}U \cap V\neq \varnothing$.
   \item \emph{point transitive}
   (PT) if there exists a point $x$ with dense orbit.
  Such a point is called \emph{transitive point}.  Notation: $x_0 \in
  Trans(X)$.
   \item \emph{densely point transitive}
   (DPT) if there exists a dense set $Y\subset X$
of transitive points.
 \end{enumerate}
\end{defin}


Of course always (DPT) implies (PT).
In general, (TT) and (PT) are independent properties.
For a detailed discussion of transitivity conditions (for
cascades) see a review paper by Kolyada and Snoha \cite{Ko}.

%



As usual, $X$ is \emph{perfect} means that $X$ is a space without
isolated points. Assertions (1) and (2) in the following
proposition are very close to Silverman's observation
\cite[Proposition 1.1]{Sil} (for cascades).

\begin{prop} \label{p:TTimplies}
\ben
\item
If $X$ is a perfect topological space and $S$ is an $F$-semigroup,
then (PT) implies (TT).
\item If $X$ is a Polish space then every (TT) system $(S,X)$ is (DPT) (and hence also (PT)).
\item  Every (DPT) system $(S,X)$ is (TT).
\een
\end{prop}
\begin{proof} (1)
Let $x$ be a transitive point with orbit $Sx$. Now, let $U$ and
$V$ be nonempty open subsets of $X$. There exists $s_1 \in S$ such
that $s_1x \in V.$  The subset $S \setminus Ss_1$ is finite
because $S$ is an almost F-group. Since $X$ is perfect, removing
the finite subset $(S \setminus Ss_1)x$ from the dense subset $Sx$
we get again a dense subset. Therefore, $Ss_1x$ is a dense subset
of $X$. Then there exists $s_2\in S$ such that $s_2s_1x\in U$.
Thus $s^{-1}_2U\cap V\neq \varnothing.$ By Definition \ref{Top
trans}.1 this means that $(S,X)$ is a (TT) dynamical system.

(2) If $(S,X)$ is topologically transitive, then $S^{-1}U$ is a
dense subset of $X$ for every open set $U$. We know that $X$ is
Polish. Then there exists a countable open base $\mathcal{B}$ of
the given topology. By the Baire theorem, $\bigcap\{S^{-1}U: \ U
\in \mathcal{B} \}$ is dense in $X$ and every point of this set is
a transitive point of the dynamical system $X$.

(3) Let $U$ and $V$ be nonempty open subsets in $X$. Since the set
$Y$ of point transitive points is dense in $X$, it intersects $V$.
Therefore, we can choose a transitive point $y\in V.$ Now by the
transitivity of $y$ there exists $s\in S$ such that $sy$ belongs
to $U$. Hence, $sy$ is a common point of $U$ and $sV$.
\end{proof}

%


%
%
%

\begin{lemma}  \label{lemma}
Let $(X,d)$ be a metric $S$-system which is (TT). Then $Eq(X)
\subset Trans(X)$.
\end{lemma}
\begin{proof}
 Let $x_0 \in Eq(X)$ and $y \in X$  .
 We have to show that
the orbit $Sx_0$ intersects the $\eps$-neighborhood
$B_{\eps}(y):=\{x \in X: d(x,y) < \eps \}$ of $y$ for every given
$\eps
>0$. Since $x_0 \in Eq(X)$ there exists a neighborhood $U$ of
$x_0$ such that $d(sx_0,sx) < \frac{\eps}{2}$ for every $(s,x) \in
S \times U$. Since $X$ is (TT) we can choose $s_0 \in S$ such that
$s_0U \cap B_{\frac{\eps}{2}}(y) \neq \emptyset$. This means that
$d(s_0x,y) < \frac{\eps}{2}$ for some $x \in U$. Then $d(s_0x_0,y)
< \eps$.
\end{proof}

\sk
\section{Minimality conditions}
\sk

The following definitions are standard for compact $X$.

\begin{defin} \label{d:minimality}
Let $X$ be a not necessarily compact S-dynamical system. \ben
\item
$X$ is called \emph{minimal}, if $\overline{Sx}=X$ for every $x\in
X$. In other words, all points of $X$ are transitive points.
\item
A point $x$ is called \emph{minimal} if the subsystem
$\overline{Sx}$ is minimal.
\item
A point $x$ is called \emph{almost periodic} if the subsystem
$\overline{Sx}$ is minimal and compact.
\item
If the set of almost periodic points is dense in $X$, we say that
$(S,X)$ satisfies the \emph{Bronstein condition}. If, in addition,
the system $(S,X)$ is (TT), we say that it is an
\emph{$M$-system.}
\item
A point $x\in X$ is a \emph{periodic point}, if $Sx$ is finite. If
(S,X) is a (TT) dynamical system and the set of periodic points is
dense in $X$, then we say that it is a \emph{$P$-system},
\cite{GW}.
\een
\end{defin}

If $X$ is compact then a point in $X$ is minimal iff it is almost
periodic.
Every periodic point is of course almost periodic. Therefore it is
also obvious that every $P$-system is an $M$-system.

 \sk For a system $(S,X)$ and a subset $B\subset X$, we use the
 following notation  $$N(x,B)=\{s\in S: sx\in B\}.$$
 The following definition is also standard.
 \begin{defin}
A subset $P\subset S$ is (left) \emph{syndetic}, if there exists a
finite set $F\subset S$ such that $F^{-1}P=S$.
 \end{defin}

 The following lemma is a slightly generalized
 version of a well known criteria for almost periodic points
 (cf. Definition \ref{d:minimality}.3)
 in compact dynamical systems. In particular,
 it is valid for every semigroup $S$.

 \begin{lemma} \label{Gott}
Let $(S,X)$ be a (not necessarily compact) dynamical system and
$x_0\in X$. Consider the following conditions:
\begin{enumerate}
\item $x_0$ is an almost periodic point.
\item For every open neighborhood $V$ of $x_0$ in $X$ there
exists a finite set $F \subset S$ such that $F^{-1}V\supseteq
Y:=\overline{Sx_0}.$
\item For every neighborhood $V$ of $x_0$ in $X$ the set $N(x_0,V)$ is
syndetic.
\item $x_0$ is a minimal point (i.e., the subsystem
$\overline{Sx_0}$ is minimal).
\end{enumerate}

Then (1) $\Rightarrow$(2) $\Rightarrow$ (3) $\Rightarrow$ (4).

If $X$ is compact then all four conditions are equivalent.

\end{lemma}
\begin{proof}
$(1) \Rightarrow (2):$ \ Suppose that $(Y,S)$ is minimal and
compact. Then for every open neighborhood $V$ of $x_0$ in $X$ and
for every $y\in Y$ there exists $s\in S$ such that $sy\in V$.
Equivalently, $y\in s^{-1}V.$ Therefore, $\bigcup_{s\in
S}s^{-1}V\supseteq Y.$ By compactness of $Y$ we can choose a
finite set $F\subseteq S$ such that $F^{-1}V\supseteq Y.$

 $(2) \Rightarrow (3):$ \ It suffices to show that $F^{-1}N(x_0,V)=S$,
 where $F$ is a subset of $S$ defined in (2).
 Assume otherwise, so that there exists $s\in S$ such that $s\notin
 F^{-1}N(x_0,V).$ Then $sx_0\notin
 F^{-1}V$. On the other hand clearly, $sx_0\in Y$,
 contrary to our condition that $F^{-1}V\supseteq Y.$

  $(3)\Rightarrow (4):$ \ $Y=\overline{Sx_0}$ is non-empty, closed and
 invariant. It remains to show that if $y\in Y$ then $x_0\in
 \overline{Sy}.$ Assume otherwise, so that $x_0\notin \overline{Sy}.$ Choose an
 open neighborhood $V$ of $x_0$ in $X$, such that $\overline{V}\cap
 \overline{Sy}=\varnothing.$ By our assumption the set $N(x_0,V)$ is
syndetic. Therefore there is a finite set $F:=\{s_1,...,s_n\}$
 so that for each $s\in S$ some $s_isx_0\in V.$ That is each
 $sx_0$ belongs to $F^{-1}V=\bigcup_{i=1}^{n}s_i^{-1}V$ for every $s \in S$.
 Hence, $Sx_0\subseteq
 \bigcup_{i=1}^{n}s_i^{-1}V.$ Then
 $$
 y \in \overline{Sx_0} \subset \overline{\bigcup_{i=1}^{n}s_i^{-1}V} = \bigcup_{i=1}^{n} \overline{s_i^{-1}V} \subset
 \bigcup_{i=1}^{n}s_i^{-1}\overline{V}.
 $$
 But then $Sy\cap \overline{V}\neq
 \varnothing$ contrary to our assumption.

If $X$ is compact then by Definition \ref{d:minimality} it follows
that (4) $\Rightarrow$ (1).
\end{proof}

\sk
\section{Sensitivity and other conditions}
\sk

\begin{prop} \label{equic-p-gener1}
 Let $S$ be an $C$-semigroup.
Assume that $(X,d)$ is a point transitive (PT) $S$-system such
that $Eq(X) \neq \emptyset$.
 Then every
transitive point is an equicontinuity point. That is, $Trans(X)
\subset Eq(X)$ holds.
\end{prop}
\begin{proof}
Let $y$ be a transitive point and $x\in Eq(X)$ be an
equicontinuity point. We have to show that $y \in Eq(X)$. For a
given $\eps >0$ there exists a neighborhood $O(x)$ of $x$ such
that
$$
d(sx'',sx')<\epsilon \quad \quad \forall \ s\in S \hsk \hsk
\forall \ x',x''\in O(x).
$$
Since $y$ is a transitive point then there exists $s_{0}\in S$
such that $s_{0}y\in O(x).$ Then $O(y):=s_{0}^{-1}O(x)$ is a
neighborhood of $y$. We have
$$
d(ss_{0}y',ss_{0}y'')<\epsilon \quad \quad \forall \ s\in S \hsk
\hsk \forall \ y', y''\in O(y).
$$

Since $S$ is a $C$-semigroup the subset $M:=\overline{S\backslash
Ss_{0}}$ is compact. Hence by Lemma \ref{com-of-S} it acts
equicontinuously on $X$. We can choose a neighborhood $U(y)$ of
$y$ such that
$$
d(ty',ty'')<\epsilon \quad \quad \forall \ t\in M \hsk \hsk
\forall \ y', y''\in U(y).
$$
Then $V:=O(y)\cap U(y)$ is a neighborhood of $y$. Since $S=M \cup
Ss_0$ we obtain that $ d(sy',sy'')<\epsilon$ for every $s\in S$
and $y',y''\in V.$ This proves that $y \in Eq(X)$.
\end{proof}

\begin{prop} \label{equic-p-gener2} Let $S$ be an
$C$-semigroup. Assume that $(X,d)$ is a metric $S$-system which is
minimal and $Eq(X) \neq \emptyset$.
Then $X$ is equicontinuous.
\end{prop}
\begin{proof}
If $(S,X)$ is a minimal system then $Trans(X)=X$. Then if $Eq(X)
\neq \emptyset$ every point is an equicontinuity point by
Proposition \ref{equic-p-gener1}. Thus, $Eq(X)=X$.
\end{proof}

\begin{prop} \label{eqnotempty}
Let $S$ be an $C$-semigroup. Assume that $(X,d)$ is a Polish (TT)
$S$-system. Then $X$ is almost equicontinuous if and only if
$Eq(X) \neq \emptyset$.
\end{prop}
\begin{proof} $X$ is (DPT) by Proposition \ref{p:TTimplies}.2.
That is, $Trans(X)$ is dense in $X$. Assuming that $Eq(X) \neq
\emptyset$ we obtain by Proposition \ref{equic-p-gener1} that
$Trans(X) \subset Eq(X)$. It follows that $Eq(X)$ is also dense in
$X$. Thus, $X$ is almost equicontinuous. This proves "if" part.
The remaining direction is trivial.
\end{proof}

\sk

The following natural definition plays a fundamental role in many
investigations about chaotic systems. The present form is a
generalized version of existing definitions for cascades (see also
\cite{GW, GM}).

\begin{defin} (sensitive dependence on initial conditions)
A metric $S$-system $(X,d)$ is \emph{sensitive} if it satisfies
the following condition: there exists a (\emph{sensitivity
constant}) $c > 0$ such that for all $x\in X$ and all $\delta>0$
there are some $y\in B_{\delta}(x)$ and $s\in S$ with $d(sx,sy)>
c$.

We say that $(S,X)$ is \emph{non-sensitive} otherwise.
\end{defin}

\begin{prop} \label{AE=NS-gen}
Let $S$ be an $C$-semigroup. Assume that $(X,d)$ is a (TT) Polish
$S$-system. Then the system is almost equicontinuous if and only
if it is non-sensitive.
\end{prop}
\begin{proof}
Clearly an almost equicontinuous system is always non-sensitive.

Conversely, the non-sensitivity means that for every $n \in \N$
there exists a nonempty open subset $V_n \subset X$ such that
$$
diam(sV_n) <
\frac{1}{n}  \quad \forall \ (s,n) \in S \times \N.
$$
Define $$U_{n}:=S^{-1}V_n  \quad \quad R:=\bigcap_{n\in
\N}U_{n}.$$ Then every $U_n$ is open. Moreover, since
 $X$ is (TT),
for every nonempty open subset $O \subset X$ there exists $s \in
S$ such that $O \cap s^{-1}U_n \neq \emptyset$. This means that
every $U_n$ is dense in $X$.
 Consequently, by Baire theorem (making use that $X$ is Polish),
 $R$ is also dense. It is enough
 now to show that $R \subset Eq(X)$.
 Suppose $x\in R$ and $\epsilon>0$. Choose $n$ so that
$\frac{1}{n} < \epsilon$, then $x\in U_{n}$ implies the existence
of $s_{0}\in S$ such that $s_{0}x\in V_n.$ Put
$V=s_{0}^{-1}V_{\frac{1 }{n}}.$ Therefore for $y\in V$ and every
$s:=s's_0 \in Ss_{0}$ we get
$$
d(sx,sy)=d(s's_{0}x,s's_{0}y)<\frac{1}{n}<\epsilon.
$$
But $S\setminus Ss_{0}$ is relatively precompact set in $S$
because $S$ is a $C$-semigroup. Then by Lemma \ref{com-of-S} the
set $S\setminus Ss_{0}$ acts on $(X,d)$ equicontinuously. We have
an open neighborhood $O$ of $x$ such that for all $y\in O$ and for
every $s\in \overline{S\setminus Ss_{0}}$ holds $d(sx,sy)<
\epsilon.$ Define an open neighborhood $M:=O\cap V$ of $x$. Then
$d(sx,sy)<\epsilon$ for every $s\in S$ and all $y\in M.$ Thus, $x
\in Eq(X)$.
\end{proof}

\begin{thm}\label{Thm}
Let $(X,d)$ be a 
Polish $S$-system where $S$ is a C-semigroup. If $X$ is
an $M$-system and $Eq(X) \neq \emptyset$ then $X$ is minimal and
equicontinuous.
\end{thm}
\begin{proof}
Let $x_{0}\in X$ be an equicontinuity point. Since every M-system
is (TT), by Lemma \ref{lemma} we know that $x_0 \in Trans(X)$.
Thus, $\overline{Sx_0}=X$. Therefore, for the minimality of $X$ it
is enough to show that $x_0$ is a minimal point.

Since $x_0 \in Eq(X)$, given $\epsilon>0,$ there exists $\delta
>0$ such that $0<\delta<\frac{\epsilon}{2}$ and $x\in
B_{\delta}(x_{0})$ implies $d(sx_{0},sx)< \frac{\epsilon}{2}$ for
every $s\in S.$ Since $X$ is an M-system the set $Y$ of all almost
periodic points is dense. Choose $y \in B_{\delta}(x_{0}) \cap Y$.
Then the set $$N(y,B_{\delta}(x_{0})): =\{s\in S: sy \in
B_{\delta}(x_{0})\}$$ is a syndetic subset of $S$ by Lemma
\ref{Gott}. Clearly, $N(y,B_{\delta}(x_{0}))$ is a subset of the
set
$$N(x_{0},B_{\epsilon}(x_{0}))=\{s\in S:d(sx_{0},x_{0}) \leq
\eps\}.$$ Then $N_{\eps}:=N(x_{0},B_{\epsilon}(x_{0}))$ is also
syndetic (for every given $\eps >0$). Using one more time Lemma
\ref{Gott} we conclude that $x_{0}$ is a minimal point, as
desired. Now the equicontinuity of $X$ follows by Proposition
\ref{equic-p-gener2}.
\end{proof}

\begin{thm} \label{main}
Let $(X,d)$ be a Polish $S$-system where $S$ is a C-semigroup. If
$X$ is an $M$-system which is not minimal or not equicontinuous.
Then $X$ is sensitive.
\end{thm}
\begin{proof}
 If $X$ is non-sensitive then by Proposition
\ref{AE=NS-gen} the system is almost equicontinuous. Theorem
\ref{Thm} implies that $X$ is minimal and equicontinuous. This
contradicts our assumption.
\end{proof}

\sk

Now if the action is a cascade $(X,f)$ or if $S$ is a topological
group (both are the case of C-semigroups, see Example \ref{e:C})
then we get, as a direct corollary, assertions (2) and (3) of
Theorem \ref{facts}. The assertion (1) is also covered in the case
of Polish phase spaces $X$. Furthermore the main results are valid
for a quite large class of actions including the actions of
one-parameter semigroups on Polish spaces.



\sk


\bibliographystyle{amsplain}

\end{document}